\documentclass[11pt]{amsart}
\usepackage{amsmath,amsthm, amscd, amssymb, amsfonts}

\usepackage{hhline}

\newcommand{\ord}{\operatorname{ord}}
\newcommand{\oper}{\#}

\newcommand\toba{{\mathfrak B }}
\newcommand\trasp{\pi}

\newcommand{\trid}{\triangleright}

\newcommand{\ku}{\mathbb C}

\newcommand{\Z}{{\mathbb Z}}
\newcommand{\N}{{\mathbb N}}

\newcommand{\C}{{\mathcal C}}
\newcommand{\D}{{\mathcal D}}

\newcommand{\Oc}{{\mathcal O}}

\newcommand\sgn{\operatorname{sgn}}

\theoremstyle{plain}

\newtheorem{maintheorem}{Theorem}
\newtheorem{lema}{Lemma}[section]
\newtheorem{theorem}[lema]{Theorem}

\newtheorem{prop}[lema]{Proposition}

\theoremstyle{definition}

\theoremstyle{remark}
\newtheorem{obs}[lema]{Remark}

\newcommand\st{\mathbb S_3}
\newcommand\sk{\mathbb S_4}
\newcommand\sn{\mathbb S_n}
\newcommand\sop{\mathbb S_p}
\newcommand\sm{\mathbb S_m}

\def\pf{\begin{proof}}
\def\epf{\end{proof}}

\theoremstyle{remark}

\begin{document}

\renewcommand{\baselinestretch}{1.2}

\thispagestyle{empty}
\title[pointed Hopf algebras with coradical $\ku\sn$]
{On pointed Hopf algebras associated to some conjugacy classes in
$\sn$}
\author[Andruskiewitsch and Zhang]{ Nicol\'as Andruskiewitsch and SHOUCHUAN ZHANG}
\address{
FaMAF, Universidad Nacional de C\'ordoba. CIEM -- CONICET
\newline \indent (5000) Ciudad Universitaria, C\'ordoba, Argentina} \email{andrus@mate.uncor.edu}
\address{
Department of Mathematics, Hunan University
\newline \indent Changsha 410082, P.R. China}\email{z9491@yahoo.com.cn}
\thanks{Results of this paper were obtained during a visit of the first author
to the Hunan University, Changsha (China).  The work of N. A. was
partially supported by CONICET, Fund. Antorchas, Agencia C\'ordoba
Ciencia, TWAS (Trieste), ANPCyT and Secyt (UNC)}

\date{\today. 2000 Mathematics Subject Classification:
      16W30}
\begin{abstract} We show that any pointed Hopf algebra with
infinitesimal braiding associated to the conjugacy class of
$\pi\in \sn$ is infinite-dimensional, if either the order of $\pi$
is odd, or all cycles in the decomposition of $\pi$ as a product
of disjoint cycles have odd order except for exactly two
transpositions.
\end{abstract}

\maketitle

\section*{Introduction}

The purpose of this article is to contribute to the classification
of finite-dimensional complex pointed Hopf algebras $H$ with $G(H)
= \sn$. Although substantial progress has been done in the
classification of finite-dimensional complex pointed Hopf algebras
with abelian group \cite{AS-05}, not much is known about the
non-abelian case. Our approach fits into the framework of the
Lifting Method \cite{AS1, AS-cambr}; we shall freely use the
notation and results from \emph{loc. cit}. Given a finite group
$G$, the key step in the classification of finite-dimensional
complex pointed Hopf algebras $H$ with $G(H) = G$ is the
determination of all Yetter-Drinfeld modules $V$ over the group
algebra of $G$ such that the Nichols algebra is
finite-dimensional. If $G$ is abelian, this reduces to the study
of Nichols algebras of diagonal type; general results were reached
in this situation in \cite{AS2, H3} for braided vector spaces of
Cartan type. If $G$ is not abelian, then very few examples have
been computed explicitly in the literature, see \cite{FK, ms, G1,
AG2}. Our viewpoint in the present paper is to deduce that some
Nichols algebras over non-abelian groups are infinite-dimensional
from those results on Nichols algebras of diagonal type. This idea
appeared first in \cite{G1}. More precisely, recall that an
irreducible Yetter-Drinfeld module over the group algebra of $G$
is determined by a conjugacy class $\C$ of $G$ and an irreducible
representation $\rho$ of the centralizer $G^s$ of a fixed $s\in
\C$. We seek for conditions on $\C$ and $\rho$ implying that the
dimension of the Nichols algebra $\toba(\C, \rho)$ is
\emph{infinite}. We concentrate on the central example $G = \sn$.
Up to now, the only pairs $(\C, \rho)$ for $\sn$ such that the
dimension $\toba(\C, \rho)$ is known are as follows.
\begin{itemize}
    \item If either $\C$ or $\rho$ is trivial, then the dimension
of $\toba(\C, \rho)$ is infinite (this is easy and holds for any
$G$).
    \item If $\C = \Oc_2$ is the class of all transpositions, $\rho = \chi_{\pm}$
    is some character
    and $n=3,4,5$, then the dimension of $\toba(\C, \rho)$ is
    finite. See \cite{FK, ms}.
    \item If $\C = \Oc_4$ is the class of all $4$-cycles in $\sk$ and $\rho= \chi_{-1}$
    is some character, then the dimension of $\toba(\C, \rho)$ is
    finite \cite[Th. 6.12]{AG2}.
\end{itemize}

It is a hard open problem to determine if the dimension of
$\toba(\Oc_2, \chi_{\pm})$ is finite or not for $n\ge 6$. The
technique in the present paper does not contribute to this
question, since all braided subspaces of diagonal type of
$M(\Oc_2, \chi_{\pm})$ give rise to Nichols algebras of finite
dimension (actually exterior algebras). We are chiefly concerned
with other representations or other orbits, hoping that a better
understanding of the special features of $M(\Oc_2, \chi_{\pm})$
(by contrast with the other orbits) would permit to isolate the
conditions for the right consideration of this hard problem.

\medbreak The main results in this paper are summarized in the
following statement.

\begin{maintheorem}\label{sthree} Let $\pi\in \sn$. If either one of
the following conditions holds

(i) the order of $\pi$ is odd, or

(ii)  all cycles in $\pi$ have odd order except for exactly two
transpositions,

then $\dim \toba(\Oc_{\pi}, \rho) = \infty$ for any $\rho\in
\widehat{\sn^{\pi}}$.
\end{maintheorem}

See Theorems \ref{nichols-oj-impar} and \ref{dosdos}. We also
apply our main result to determine all irreducible Yetter-Drinfeld
modules over $\st$ or $\sk$ whose Nichols algebra is
finite-dimensional. See Theorem \ref{nichols-sthree}.

\section{Preliminaries and conventions}\label{conventions}

Our references for pointed Hopf algebras are \cite{M},
\cite{AS-cambr}. The set of isomorphism classes of irreducible
representations of a finite group $G$ is denoted $\widehat{G}$;
thus, the group of characters of a finite abelian group $\Gamma$
is denoted $\widehat{\Gamma}$. We shall often confuse a
representant of a class in $\widehat{G}$ with the class itself. If
$V$ is a $\Gamma$-module then $V^\chi$ denotes the isotypic
component of type $\chi\in \widehat{\Gamma}$. It is useful to
write $g\trid h = ghg^{-1}$, $g, h\in G$.

\subsection{Notation on the groups $\sn$}\label{conventionss3}

Recall that the type of a permutation $\trasp \in \sn$ is a symbol
$(1^{m_1}, 2^{m_2}, \dots, n^{m_n})$ meaning that in the
decomposition of $\trasp$ as product of disjoint cycles, there are
$m_j$ cycles of length $j$, $1\le j \le n$. We may omit $j^{m_j}$
when $m_j = 0$. The conjugacy class $\Oc_\trasp$ of $\trasp$
coincides with the set of all permutations in $\sn$ with the same
type as $\trasp$; we denote it by $\Oc_{1^{m_1}, 2^{m_2}, \dots,
n^{m_n}}$ or by $\Oc_{2, \dots, 2, 3 \dots, 3, \dots}$. For
instance, we denote by $\Oc_j$ the conjugacy class of $j$-cycles
in $\sn$, $2\le j \le n$. Also, $\Oc_{2,2}$ is the conjugacy class
of $(12)(34)$ and so on. If $\trasp \in \sn$ and $n<m$ then we
also denote by $\trasp$ the natural extension to $\sm$ that fixes
all $i>n$. If some emphasis is needed, we add a superscript $n$ to
indicate that we are taking conjugacy classes in $\sn$, like
$\Oc_j^n$ for the conjugacy class of $j$-cycles in $\sn$. It is
well-known that the isotropy subgroup $\sn^{\pi}$ is isomorphic to
a product
$$\sn^{\pi} \simeq T_{1}\dots T_{n}$$ where $T_i =  \Gamma_i
\rtimes \mathbb S_{m_i}$, $1\le i \le n$. Here $\Gamma_i \simeq
(\Z/ i)^{m_i}$ is generated by the $i$-cycles in $\pi$ and
$\mathbb S_{m_i}$ permutes these cycles. Hence any $\rho \in
\widehat{\sn^{\pi}}$ is of the form $\rho_1\otimes \dots \otimes
\rho_n$ where $\rho_i \in \widehat{T_i}$.

If $X$ is a subset of $\{1, \dots, n\}$ then $\mathbb S_X$ denotes
the subgroup of bijections in $\sn$ that fix pointwise $\{1,
\dots, n\} - X$. We shall use the following notation for
representations of subgroups of $\sn$. Let $\omega_j$ be a fixed
primitive $j$-th root of 1, $2\le j \le n$. Let $\tau_j = (123
\dots j)$.
\begin{align*}
\varepsilon &= \text{ trivial character}; \\
\sgn &= \text{ sign character of  $\mathbb S_X$, } X\subset \sn;
\\ \chi_j &= \text{  character  of $\langle\tau_j\rangle\simeq \Z_j$ given by }
\chi(\tau_j) = \omega_j.
\end{align*}

\subsection{Yetter-Drinfeld modules over the group algebra of a finite
group}\label{conventionsyd}

A Yetter-Drinfeld module over a finite group $G$ is a left
$G$-module and left $\ku G$-comodule $M$ such that
$$
\delta(g.m) = ghg^{-1} \otimes g.m, \qquad m\in M_h, g, h\in G.
$$
Here $M_h = \{m\in M: \delta(m) = h\otimes m\}$; clearly, $M =
\oplus_{h\in G} M_h$. Yetter-Drinfeld modules over $G$ are
completely reducible. Also, irreducible Yetter-Drinfeld modules
over $G$ are parameterized by pairs $(\C, \rho)$ where $\C$ is a
conjugacy class and $\rho$ is an irreducible representation of the
isotropy subgroup $G^s$ of a fixed point $s\in \C$. As usual,
$\deg \rho$ is the dimension of the vector space $V$ affording
$\rho$. We denote the corresponding Yetter-Drinfeld module by
$M(\C, \rho)$; see \cite{DPR, Wi}, and also \cite{AG1}. Since
$s\in Z(G^s)$, the Schur Lemma says that
\begin{equation}\label{schur} s \text{ acts by a scalar $q_{ss}$
on } V.
\end{equation}

Here is a precise description of the Yetter-Drinfeld module $M(\C,
\rho)$. Let $t_1 = s$, \dots, $t_{M}$ be a numeration of $\C$ and
let $g_i\in G$ such that $g_i \trid s = t_i$ for all $1\le i \le
M$. Then  $M(\C, \rho) = \oplus_{1\le i \le M}g_i\otimes V$. Let
$g_iv := g_i\otimes v \in M(\C,\rho)$, $1\le i \le M$, $v\in V$.
If $v\in V$ and $1\le i \le M$, then the action of $g\in G$ and
the coaction are given by
$$
\delta(g_iv) = t_i\otimes g_iv, \qquad g\cdot (g_iv) =
g_j(\gamma\cdot v),
$$
where $gg_i = g_j\gamma$, for some $1\le j \le M$ and $\gamma\in
G^s$. The explicit formula for the braiding is then given by
\begin{equation} \label{yd-braiding}
c(g_iv\otimes g_jw) = t_i\cdot(g_jw)\otimes g_iv = g_h(\gamma\cdot
v)\otimes g_iv\end{equation} for any $1\le i,j\le M$, $v,w\in V$,
where $t_ig_j = g_h\gamma$ for unique $h$, $1\le h \le M$ and
$\gamma \in G^s$.

\subsection{On Nichols algebras}\label{conventionsnichols}

Let $(V, c)$ be a braided vector space, that is $V$ is a vector
space and $c: V\otimes V\to V\otimes V$ is a linear isomorphism
satisfying the braid equation. Then $\toba(V)$ denotes the Nichols
algebra of $V$, see the precise definition in \cite{AS-cambr}. Let
$G$ be a finite group, $\C$ a conjugacy class, $s\in \C$ and
$\rho\in \widehat{G^s}$. The Nichols algebra of $M(\C, \rho)$ will
be denoted simply by $\toba(\C, \rho)$. We collect some general
facts on Nichols algebras for further reference.

\begin{obs}\label{trivialbraiding}
Let $(V,c)$ be a braided vector space. If $W$ is a subspace of $V$
such that $c(W\otimes W) = W\otimes W$ then $\toba(W) \subset
\toba(V)$. Thus $\dim \toba(W) =\infty \implies \dim \toba(V)
=\infty$. In particular, if there exists $v\in V - 0$ such that
$c(v\otimes v ) = v \otimes v$ then $\dim \toba(V) =\infty$. This
is the case in $\toba(\C, \rho)$ when either the orbit $\C$ is
trivial or the representation $\rho$ is trivial, or even if the
scalar $q_{ss}$ defined in \eqref{schur} is 1.
\end{obs}

A braided vector space $(V,c)$ is of \emph{diagonal type} if there
exists a basis $v_1, \dots, v_{\theta}$ of $V$ and non-zero
scalars $q_{ij}$, $1\le i,j\le \theta$, such that
$$
c(v_i\otimes v_j) = q_{ij} v_j\otimes v_i, \qquad \text{for all }
1\le i,j\le \theta.
$$
A braided vector space  $(V,c)$ is of \emph{Cartan type} if it is
of diagonal type, and there exists $a_{ij} \in \Z$, $-\ord q_{ii}
< a_{ij} \leq 0$ such that
$$
q_{ij}q_{ji} = q_{ii}^{a_{ij}}
$$
for all $1\le i\neq j\le \theta$.  Set $a_{ii}=2$ for all $1\le
i\le \theta$. Then $(a_{ij})_{1\le i,j\le \theta}$ is a
generalized Cartan matrix. The following important result was
proved in \cite[Th. 4]{H3}, showing that some hypotheses in
\cite[Th. 1]{AS2} were unnecessary.

\begin{theorem}\label{cartantype}
Let $(V,c)$ be a braided vector space of Cartan type. Assume that
$q_{ii}\neq 1$ is root of 1 for all $1\le i\le \theta$. Then $\dim
\toba(V) < \infty$ if and only if the Cartan matrix is of finite
type. \qed\end{theorem}

\section{On Nichols algebras over $\sn$}\label{nichols-sn}

In this section we state some general results about Nichols
algebras over $\sn$. Our main idea is to find out suitable braided
subspaces of a Yetter-Drinfeld module over $\ku \sn$ that are
diagonal of Cartan type, so that Theorem \ref{cartantype} applies.
In what follows, $G$ is a finite group, $s\in G$, $\C$ is the
conjugacy class of $s$, $\rho\in \widehat{G^s}$, $\rho: G^s \to
GL(V)$. Recall the scalar $q_{ss}$ defined in \eqref{schur}.

\begin{prop}\label{matias} \cite[Lemma 3.1]{G1}
Assume that $\dim\toba(\C, \rho)< \infty$. Then
\begin{itemize}
    \item $\deg \rho > 2$ implies $q_{ss} = -1$.
    \item $\deg \rho = 2$ implies $q_{ss} = -1$, $\omega_3$ or $\omega_3^2$.\qed
\end{itemize}
\end{prop}

\subsection{Nichols algebras corresponding to permutations of odd order}\label{nichols-odd}

We begin by a general way of finding braided subspaces of rank 2.
Assume that
\begin{equation}\label{hypothesis-involution}
\text{there exists  $\sigma\in G$ such that }\sigma s \sigma^{-1}
= s^{-1}\neq s.
\end{equation}
Hence $s^{-1}\in \C$, \emph{i. e.} $\C$ is real. Under this
hypothesis, we prove the following result that, unlike Proposition
\ref{matias}, does not assume any restriction on $\deg \rho$.

\begin{lema}\label{odd} If $\dim\toba(\C, \rho)< \infty$ then $q_{ss} = -1$, and $s$ has even order.
\end{lema}

\pf Let $N = \ord q_{ss}$; clearly $N>1$. Let $t_1 = s$, $t_2 =
s^{-1}$, $g_1 = e$, $g_2 = \sigma$. By
\eqref{hypothesis-involution}, $\sigma s^{-1} \sigma^{-1} = s$;
thus $\sigma s^{-1} = s\sigma$ and $\sigma s = s^{-1}\sigma$.
Together with \eqref{schur} and \eqref{yd-braiding}, this shows
that for any $v,w\in V$
\begin{align*}
c(g_1 v \otimes g_1 w) &= g_1 (s\cdot w) \otimes g_1 v = q_{ss} \,
g_1 w \otimes g_1 v, \\ c(g_1 v \otimes g_2 w) &= g_2(s^{-1}\cdot
w) \otimes g_1 v =
q_{ss}^{-1} \, g_2 w \otimes g_1 v,\\
c(g_2 v \otimes g_1 w) &= g_1 (s^{-1}\cdot w)  \otimes g_2 v=
q_{ss}^{-1} \, g_1 w \otimes g_2 v,
\\ c(g_2 v \otimes g_2 w) &= g_2 (s\cdot w) \otimes g_2v
= q_{ss} \, g_2 w \otimes g_2 v.
\end{align*}
Let $v\in V$, $v\neq 0$. Then the subspace of $M(\C,\rho)$ spanned
by $v_{1} := g_1v$, $v_{2} := g_2v$ is a braided subspace of
Cartan type: $c(v_{i} \otimes v_{j}) = q_{ij} v_{j}\otimes v_{i}$,
for $1\le i,j\le 2$, where $\begin{pmatrix} q_{11} & q_{12} \\
q_{21} & q_{22}  \end{pmatrix} =  \begin{pmatrix} q_{ss} & q_{ss}^{-1} \\
q_{ss}^{-1} & q_{ss} \end{pmatrix}$,
with Cartan matrix $\begin{pmatrix} 2 & a_{12} \\
a_{21} & 2 \end{pmatrix}$, $a_{12}= a_{21} \equiv -2 \mod N$. Now
$a_{12} = a_{21} =0$ or $-1$, by hypothesis and Theorem
\ref{cartantype}. Thus, necessarily $a_{12} = a_{21} =0$ and $N =
2$. \epf

\begin{lema}\label{nichols-cycle} Let $\trasp\in \sn$, $\ord \pi > 2$, $\rho\in
\widehat{\sn^\trasp}$. If $\dim\toba(\Oc_\pi,\rho) < \infty$ then
$q_{\trasp, \trasp} = -1$.
\end{lema}

\pf It is well-known that \eqref{hypothesis-involution} holds for
any $\trasp\in \sn$.  Namely, assume that $\trasp = t_j$ for some
$j$ and take
$$g_2 =
  \begin{cases} &(1 \, j-1) (2\, j-2) \cdots (k-1\, k+1), \text{ if } j = 2k \text{ is even}, \\
 &(1 \, j-1) (2\, j-2) \cdots (k\, k+1), \,\,\,\quad \,\text{ if } j = 2k + 1 \text{ is odd.}
  \end{cases}$$
It is easy to see that  $g_2t_jg_2^{-1} = t_j^{-1}$. The general
case follows using that any $\trasp$ is a product of disjoint
cycles. We conclude from Lemma \ref{odd}. \epf

\begin{theorem}\label{nichols-oj-impar} If $\trasp\in \sn$ has odd order
then $\dim\toba(\Oc_\trasp, \rho) = \infty$ for any $\rho \in
\widehat{\sn^{\trasp}}$. \qed
\end{theorem}

If $n$ es even, the Nichols algebras $\toba(\Oc_n^n, \rho)$ cannot
be treated by similar arguments as above. For, the isotropy
subgroup $\sn^{\tau_n}$ is cyclic of order $n$ and we can assume
that $\rho(\tau_n)= -1$. Assume that $\tau \in \Oc_n^n \cap
\sn^{\tau_n}$  then $\tau = \tau_n^j$ with $(n,j)=1$ hence
$\rho(\tau)= -1$. But there are Nichols algebras like these that
are finite-dimensional. For instance, the Nichols algebra
$\toba(\Oc_4, \chi_4^2)$ was computed in \cite[Th. 6.12]{AG2} and
has dimension 576.

\subsection{A reduction argument}\label{nichols-reduction}

We now discuss a general reduction argument. Let $n,p\in \N$ and
let $m= n+p$. Let $\oper: \sn \times \sop \to \sm$ be the group
homomorphism given by
$$
\pi\oper\tau (i) =
\begin{cases}
\pi(i), & 1\le i \le n; \\
\tau(i-n) +n, & n+1\le i \le m,
\end{cases}
$$
$\pi\in \sn$, $\tau\in \sop$. If also $g\in \sn$, $h\in \sop$ then
$(g\oper h) \trid (\pi\oper \tau) = (g\trid\pi)\oper (h\trid
\tau)$. Thus $\Oc_{\pi\oper \tau} \supset \Oc_{\pi}\oper
\Oc_{\tau}$. Let us say that $\pi$ and $\tau$ are
\emph{orthogonal}, denoted $\pi\perp \tau$, if there is no $j$
such that both $\pi$ and $\tau$ contain a $j$-cycle. In other
words, if $\pi$ has type $(1^{a_1}, 2^{a_2}, \dots, n^{a_n})$ and
$\tau$ has type $(1^{b_1}, 2^{b_2}, \dots, p^{b_p})$ then either
$a_j = 0$ or $b_j = 0$ for any $j$. If $\pi\perp \tau$ then
$\sm^{\pi\oper\tau} = \sn^{\pi}\oper\sop^{\tau}$, say by a
counting argument. Hence any $\mu \in
\widehat{\sm^{\pi\oper\tau}}$ is of the form $\mu = \rho \otimes
\lambda$, for unique $\rho \in \widehat{\sn^{\trasp}}$, $\lambda
\in \widehat{\sop^{\tau}}$. Say that $V$, $W$ are the vector
spaces affording $\rho$, $\lambda$. Let $q_{\tau\tau}$, resp.
$q_{\pi\pi}$, be the scalar as in \eqref{schur} for the
representation $\rho$, resp. $\lambda$.

\begin{lema}\label{nichols-descomp} Let $\pi\in \sn$, $\tau\in
\sop$.

(1) Assume that $\pi\perp \tau$ and $\ord(\pi\oper \tau) > 2$. Let
$\mu \in \widehat{\sm^{\pi\oper\tau}}$ of the form $\mu = \rho
\otimes \lambda$, for $\rho \in \widehat{\sn^{\trasp}}$, $\lambda
\in \widehat{\sop^{\tau}}$.  If $\dim\toba(\Oc_{\pi\oper\tau},\mu)
< \infty$ then $q_{\pi\pi}q_{\tau\tau} = -1$.

(2) Assume that the orders of $\pi$ and $\tau$ are relatively
prime, that $\ord q_{\tau\tau}$ is odd and that $\ord(\pi\oper
\tau) > 2$. If $\dim\toba(\Oc_{\pi\oper\tau},\mu) < \infty$ then
$q_{\tau\tau} = 1$.
\end{lema}

Clearly, if the orders of $\pi$ and $\tau$ are relatively prime
then $\pi\perp \tau$. On the other hand, if $\ord(\pi\oper \tau) =
2$ and the orders of $\pi$ and $\tau$ are relatively prime, then
we can assume $\tau = e$, and $q_{\tau\tau} = 1$ anyway.

\pf Since $q_{\pi\oper \tau, \pi\oper\tau} = q_{\pi,
\pi}q_{\tau\tau}$, (1) follows from Lemma \ref{nichols-cycle}.
Then (2) follows from (1).  \epf

Our aim is to obtain information on $\toba(\Oc_{\pi\oper \tau},
\mu)$ from $\toba(\Oc_{\pi}, \rho)$ and $\toba(\Oc_{\tau},
\lambda)$. For this, we fix $g_1 = e, \dots, g_{P} \in \sn$, $h_1
= e, \dots, h_{T}\in \sop$,  such that
\begin{align*}
g_1 \trid \pi &= \pi, \dots, g_{P}\trid \pi \text{ is a numeration
of } \Oc_{\pi},
\\ h_1 \trid \tau &= \tau, \dots, h_{T}\trid \tau \text{ is a numeration  of }
\Oc_{\tau}.
\end{align*}

Then we can extend $(g_i\oper h_j) \trid (\pi\oper \tau)$, $1\le i
\le P$, $1\le j \le T$ to a numeration of $\Oc_{\pi\oper \tau}$.
Let $v,u\in V$, $w,z\in W$, $1\le i,k\le P$, $1\le j,l \le T$.
Then the braiding in $M(\Oc_{\pi\oper \tau},\rho \otimes \lambda)$
has the form
\begin{multline} \label{yd-braiding-oper}
c\left((g_i\oper h_j)(v\otimes w) \otimes (g_k\oper h_l)(u\otimes
z)\right) \\
= \left( (g_r\oper h_p)(\gamma\cdot u\otimes \beta\cdot z) \otimes
(g_i\oper h_j)(v\otimes w)\right)\end{multline} where
$(g_i\trid\pi) g_k = g_r\gamma$, $(h_j\trid\tau) h_l = h_p \beta$
for unique  $1\le r \le P$, $1\le p \le T$, $\gamma \in \sn^\pi$
and $\beta \in \sop^\tau$. Assume in \eqref{yd-braiding-oper} that
$j = l = 1$; then $p=1$, $\beta = \tau$ and
\eqref{yd-braiding-oper} takes the form
\begin{multline} \label{yd-braiding-oper2}
c\left((g_i\oper h_1)(v\otimes w) \otimes (g_k\oper h_1)(u\otimes
z)\right) \\
= q_{\tau\tau} \left( (g_r\oper h_1)(\gamma\cdot u\otimes z)
\otimes (g_i\oper h_1)(v\otimes w)\right).\end{multline}

Our aim is to spell out some consequences of formula
\eqref{yd-braiding-oper2}.

\begin{prop}\label{reduction} Let $\pi\in \sn$, $\tau\in
\sop$. Assume that the orders of $\pi$ and $\tau$ are relatively
prime and that $\ord q_{\tau\tau}$ is odd. Let $\mu \in
\widehat{\sm^{\pi\oper\tau}}$ of the form $\mu = \rho \otimes
\lambda$, for $\rho \in \widehat{\sn^{\trasp}}$, $\lambda \in
\widehat{\sop^{\tau}}$.

(1)  If $\dim\toba(\Oc_{\pi\oper\tau}, \mu) < \infty$, then
$\dim\toba(\Oc_{\pi}, \rho) < \infty$.

(2) If $\dim\toba(\Oc_{\pi}, \rho) = \infty$ for any $\rho \in
\widehat{\sn^{\pi}}$, then $\dim\toba(\Oc_{\pi\oper\tau}, \mu) =
\infty$ for any $\mu \in \widehat{\sm^{\pi\oper\tau}}$.
\end{prop}

In particular, let $\trasp\in \sn$ with no fixed points. If
$\dim\toba(\Oc^n_\trasp, \rho) = \infty$ for any $\rho \in
\widehat{\sn^{\trasp}}$, then $\dim\toba(\Oc^m_\trasp, \rho') =
\infty$ for any $\rho' \in \widehat{\sm^{\trasp}}$, $m
> n$.

\pf Assume that $\dim\toba(\Oc_{\pi\oper\tau}, \mu) < \infty$. By
Lemma \ref{nichols-descomp}, $q_{\tau\tau}=1$. Let $0 \neq w=z\in
W$. The linear map $\psi: M(\Oc_{\pi}, \rho)  \to
M(\Oc_{\pi\oper\tau}, \mu)$ given by $\psi(g_iv) = (g_i\oper
h_1)(v\otimes w)$ is a morphism of braided vector spaces because
of \eqref{yd-braiding-oper2}. Now apply Remark
\ref{trivialbraiding}. \epf

\subsection{Nichols algebras of orbits with exactly two transpositions}\label{nichols-four}

\begin{theorem}\label{dosdos} Let $\trasp\in \sn$. If $\trasp$ has type
$(1^a, 2^2, h_1^{m_1}, \dots,h_r^{m_r})$ where $h_1, \dots,h_r$
are odd, then $\dim\toba(\Oc^n_\trasp, \rho) = \infty$ for any
$\rho \in \widehat{\sn^{\trasp}}$.
\end{theorem}

For instance, if $n\geq 4$ then $\dim\toba(\Oc^n_{2,2}, \rho) =
\infty$ for any $\rho$.

\pf By Proposition \ref{reduction}, we can assume that $n=4$ and
$\pi$ is of type $(2,2)$. Let us consider the irreducible
Yetter-Drinfeld modules corresponding to $\Oc_{2,2} = \{a
=(13)(24), b = (12)(34), d = (14)(23)\}$. The isotropy subgroup of
$a$ is $\sk^a = \langle (1234), (13)\rangle \simeq\D_4$. Let $A=
(1234)$, $B = (13)$; $\Oc_{2,2} \subseteq \sk^a$ and $a = A^2$,
$b=BA$, $d= BA^3$. Hence the irreducible representations of
$\sk^a$ are (1) the characters given by $A\mapsto \varepsilon_1$,
$B\mapsto \varepsilon_2$ where $\varepsilon_j \in \{\pm 1\}$,
$1\le j \le 2$, and (2) the 2-dimensional representation
$\rho:\mathbb D_4 \to GL(2,\mathbb C) $ given by
$$\rho(A) = \begin{pmatrix} 0& -1\\ 1& 0 \end{pmatrix}, \qquad \rho(B)
= \begin{pmatrix} -1& 0\\ 0& 1 \end{pmatrix}.$$

\medbreak Let $\mu$ be a one-dimensional representation of
$\sk^a$. Then $M(\Oc_{2,2}, \mu)$ has a basis $v_a, v_b, v_d$ with
$\delta(v_a) = a \otimes v_a$, etc., and $c(v_a \otimes v_a) =
\mu(a) v_a \otimes v_a = v_a \otimes v_a$. Hence $\dim
\toba(\Oc_{2,2}, \mu) = \infty$.

\medbreak Let $\sigma_0 = 1$, $\sigma_1 = (12)$, $\sigma_2 =
(23)$. Then
$$
\sigma_1\trid a = d, \quad \sigma_2\trid a = b, \quad
\sigma_1\trid b = b, \quad \sigma_2\trid d = d.
$$
Let us consider the Yetter-Drinfeld module $M(\Oc_{2,2}, \rho)$.
We have
$$\rho(a) = \begin{pmatrix}  -1 & 0\\ 0 & -1 \end{pmatrix}, \quad \rho(b)
= \begin{pmatrix} 0& 1\\ 1& 0 \end{pmatrix}, \quad \rho(d) =
\begin{pmatrix} 0 &-1\\ -1& 0 \end{pmatrix}.$$
Let $\sigma_jv := \sigma_j \otimes v$, $v\in V$, $0\le j \le 2$.
The coaction is given by $\delta(\sigma_j v) = \sigma_j \trid a
\otimes \sigma_j v$; we need the action of the elements $a$, $b$,
$d$, which is
\begin{align*}
a\cdot \sigma_0 v &= \sigma_0 \rho(a)(v), &\qquad a\cdot \sigma_1
v &= \sigma_1 \rho(d)(v), &\qquad a\cdot \sigma_2 v &= \sigma_2 \rho(b)(v), & \\
b\cdot \sigma_0 v &= \sigma_0 \rho(b)(v), &\qquad b\cdot \sigma_1
v &= \sigma_1 \rho(b)(v), &\qquad b\cdot \sigma_2 v &= \sigma_2 \rho(a)(v), & \\
d\cdot \sigma_0 v &= \sigma_0 \rho(d)(v), &\qquad d\cdot \sigma_1
v &= \sigma_1 \rho(a)(v), &\qquad d\cdot \sigma_2 v &= \sigma_2
\rho(d)(v).
\end{align*}

\noindent Hence the braiding is given, for all $0\le j \le 2$ and
$v,w\in V$, by
\begin{align*}
c(\sigma_j v \otimes \sigma_j w) &= (\sigma_j\trid a)\cdot
\sigma_jw \otimes \sigma_j v = \sigma_j \rho(a)(w) \otimes
\sigma_j v = -\sigma_j w \otimes \sigma_j v
\end{align*}
and
\begin{align*}
c(\sigma_0 v \otimes \sigma_1 w) &= \sigma_1 \rho(d)(w) \otimes
\sigma_0 v, \qquad c(\sigma_0 v \otimes \sigma_2 w) &= \sigma_2
\rho(b)(w) \otimes \sigma_0 v,\\
c(\sigma_1 v \otimes \sigma_0 w) &= \sigma_0 \rho(d)(w) \otimes
\sigma_1 v, \qquad c(\sigma_1 v \otimes \sigma_2 w) &= \sigma_2
\rho(d)(w) \otimes \sigma_1 v, \\ c(\sigma_2 v \otimes \sigma_0 w)
&= \sigma_0 \rho(b)(w) \otimes \sigma_2 v, \qquad c(\sigma_2 v
\otimes \sigma_1 w) &= \sigma_1 \rho(b)(w) \otimes \sigma_2
v.\end{align*}

Let $v_1 = \begin{pmatrix} 1 \\ 1\end{pmatrix}$, $v_2 =
\begin{pmatrix}  1 \\ -1\end{pmatrix}$. Then $\rho(b)(v_1) = v_1$,
$\rho(b)(v_2) = -v_2$, $\rho(d)(v_1) = -v_1$, $\rho(d)(v_2) =
v_2$. Hence the braiding is diagonal of Cartan type in the basis
$$w_1 = \sigma_0v_1,\,\, w_2 = \sigma_0v_2, \,\, w_3 =
\sigma_1v_1,\,\, w_4 = \sigma_1v_2,\,\,  w_5 = \sigma_2v_1,\,\,
w_6 = \sigma_2v_2.$$

The corresponding Dynkin diagram is not connected; its connected
components are $\{1, 4, 6\}$ and  $\{2, 3, 5\}$, each of them
supporting the affine Dynkin diagram $A_2^{(1)}$. Then $\dim
\toba(\Oc_{2,2}, \rho) = \infty$ by Theorem \ref{cartantype}. \epf

\subsection{Nichols algebras over $\st$ and $\sk$}\label{do-gpd}

We  apply the main result of this paper to classify
finite-dimensional Nichols algebras over $\st$ or $\sk$ with
irreducible module of primitive elements.

\begin{table}[b]
\begin{center}
\begin{tabular}{|p{1,5cm}|p{2cm}|p{2,9cm}|p{1,7cm}|p{2,2cm}|}
\hline {\bf Orbit} &    {\bf Isotropy \newline group}& {\bf
Representation}
 & $\dim \toba(V)$ & {\bf Reference}

\\ \hline  $e$   & $\st$  &   any     & $\infty$ & Remark \ref{trivialbraiding}

\\ \hline  $\Oc_3$   & $\Z_3$  &   any     & $\infty$ &
Theorem \ref{nichols-oj-impar}

\\ \hline  $\Oc_2$   & $\Z_2$  &   $\varepsilon$     & $\infty$
& Remark \ref{trivialbraiding}

\\ \hline  $\Oc_2$   & $\Z_2$  &   $\sgn$
 & $12$ & \cite{ms}

\\ \hline
\end{tabular}
\end{center}

\caption{Nichols algebras of irreducible Yetter-Drinfeld modules
over $\st$ }\label{tablauno}

\end{table}

\begin{table}[b]
\begin{center}
\begin{tabular}{|p{1,5cm}|p{2cm}|p{2,9cm}|p{1,7cm}|p{2,2cm}|}
\hline {\bf Orbit} &    {\bf Isotropy \newline group}& {\bf
Representation}
 & $\dim \toba(V)$ & {\bf Reference}

\\ \hline  $e$ & $\sk$   & any & $\infty$& Remark \ref{trivialbraiding}

\\ \hline  $\Oc_{2,2}$ & $\D_4$ & any &$\infty$ & Theorem \ref{dosdos}

\\ \hline  $\Oc_4$ &$\Z_4 $
 &
$\varepsilon$ & $\infty$ & Remark \ref{trivialbraiding}

\\ \hline  $\Oc_4$ &$\Z_4 $
 &
 $\chi_4$ or $\chi_4^3$ & $\infty$ & Lemma
\ref{nichols-cycle}

\\ \hline   $\Oc_4$ &$\Z_4 $ & $\chi_4^2$ &  576
& \cite[6.12]{AG2}

\\ \hline  $\Oc_3$ & $\Z_3$
 & any &$\infty$ & Theorem \ref{nichols-oj-impar}

\\ \hline   $\Oc_2$ & $\Z_2\oplus \Z_2$
& $\varepsilon$ or  $\varepsilon \oplus \sgn$ &$\infty$
& Remark \ref{trivialbraiding}

\\ \hline   $\Oc_2$ & $\Z_2\oplus \Z_2$    &   $\sgn \oplus \varepsilon$ & 576 & \cite{FK}

\\ \hline   $\Oc_2$ & $\Z_2\oplus \Z_2$    & $\sgn \oplus \sgn$ & 576
& \cite{ms}

\\ \hline
\end{tabular}
\end{center}

\caption{Nichols algebras of irreducible Yetter-Drinfeld modules
over $\sk$}\label{tablados}

\end{table}

\begin{theorem}\label{nichols-sthree}
Let $M(\C, \lambda)$ be an irreducible Yetter-Drinfeld module over
$\sn$ such that $\toba(\C, \lambda)$ is finite-dimensional.

 (i). If $\sn =\st$ then $M(\C, \lambda) \simeq M(\Oc_2, \sgn)$.

(ii). If $\sn =\sk$ then $M(\C, \lambda)$ is isomorphic either to
$\toba(\Oc_4, \chi_4^2)$ or to $\toba(\Oc_2, \sgn\oplus
\varepsilon)$ or to $\toba(\Oc_2, \sgn\oplus \sgn)$.
\end{theorem}

\pf See tables \ref{tablauno} and \ref{tablados}. \epf

\subsection*{Acknowledgements} The first author thanks Istv\'an Heckenberger
and Mat\'\i as Gra\~na for interesting e-mail exchanges.


\end{document}